\documentstyle[11pt]{article}

\def\RR{{\rm I}\!{\rm R}}
\def\PP{{\rm I}\!{\rm P}}
\def\QQ{{\rm Q}\!\!\!\vrule height 6.4pt depth -0.4pt width 1pt \ }
\def\ZZ{{\rm Z}\!\!{\rm Z}}

\def\AA{{\it I}\!\!{\rm A}}
\def\NN{{{\rm I}\!{\rm N}}}

\def\I {\tilde{I}}
\def\f {\tilde{f}}

\def\A {{\cal A}}
\def\B {{\cal B}}
\def\N {{\cal N}}
\def\P {{\cal P}}
\def\Q {{\cal Q}}
\def\S {{\cal S}}
\def\U {{\cal U}}

\def\pos {{\rm \mbox{pos}}}

\def\conv {{\rm \mbox{conv}}}
\def\Vol {{\rm \mbox{Vol}}}
\def\vol {{\rm \mbox{vol\,}}}
\def\mod {{\rm \mbox{mod\, }}}

\newcommand{\comb}[2]{%
			  (^{#1}_{#2})}

\newcommand{\ssp}{\hspace{20mm}}
\newcommand{\msp}{\hspace{30mm}}
\newcommand{\lsp}{\hspace{40mm}}

\newtheorem{defn}{Definition}[section]
\newtheorem{lem}[defn]{Lemma}
\newtheorem{mainlem}[defn]{Main Lemma}

\newtheorem{thm}[defn]{Theorem}
\newtheorem{cor}[defn]{Corollary}

\newtheorem{exmp}[defn]{Example}

\newenvironment{undef}[1]%
		   {\vspace{3.3mm}
		   \noindent{\bf #1}\it}%
		   {\vspace{3.3mm}}

\newenvironment{proof}[1]{
  \trivlist \item[\hskip \labelsep{\it #1}]}{\hfill\mbox{$\Box$}
  \endtrivlist}

\begin{document}

% ----------- Titulo -----------------------------------------------

\noindent {\Large {\bf  A Sparse Effective Nullstellensatz 
			\footnote{Date of this version: 6.8.97}}}

\vspace{15mm}

\noindent {\large Mart\'\i n Sombra \footnote{Partially supported
	   by CONICET PID 3949/92, UBA CyT EX. 001 and Fundaci\'on 
	   Antorchas.}

\vspace{1mm}

\noindent {\small Departamento de Matem\'atica, 
Universidad  
		  Nacional de La Plata,\\
		  Calle 50 y 115, 1900 La Plata, Argentina.\\
		  {\tt E-mail: sombra@mate.unlp.edu.ar}}

\vspace{10mm}

% ----------- Abstract -----------------------------------------------

\noindent {\small {\bf Abstract.}  
We present bounds for the sparseness and for the degrees of the polynomials
in the Nullstellensatz. 
Our bounds depend mainly on the
unmixed volume of the input polynomial system. 
The degree bounds can substantially improve the known  ones when this 
polynomial system is sparse, and they are, in the worst case, simply 
exponential in terms of the number of variables and the maximum degree 
of the input polynomials.}

\vspace{4mm}


\vspace{4mm}


\vspace{8mm}

% ----------- Introduccion -----------------------------

\noindent {\bf Introduction}

\vspace{5mm}

We consider the problem of bounding the 
sparseness of the polynomials in the Nullstellensatz in
the case when the input is a sparse polynomial system. 

\vspace{3mm}

Let us denote by $k$ a field and by $\overline{k}$
its algebraic closure. 
As usual we denote by $\AA^n$
the affine space of $n$ dimensions over $\overline{k}$.  
Let be given polynomials $f_1, \ldots, f_s \in 
k[x_1, \ldots, x_n]$ without common 
zeros in $\AA^n$. 
Classical Hilbert's Nullstellensatz states then that there exist 
$g_1, \ldots, g_s \in k[x_1, \ldots, x_n]$ satisfying the 
B\'ezout equation
\[
\lsp 1= g_1 f_1 + \cdots + g_s f_s. \lsp (1)
\]

Let $d$ denote the maximum degree of the polynomials 
$f_1, \ldots, f_s$ and assume that $n\ge 2$. 
In the 
previous situation we have indeed that
there exist polynomials $g_1, \ldots, g_s $ verifying the
degree bound 
\[
\deg g_i f_i \le \max\{3, d\}^n. 
\]
This result is due to Koll\'ar, and it is optimal 
provided that $d\ge 3$ \cite{Kollar88}. 
In the 
case when $d=2$ and the characteristic of $k$ is different from two, 
Sabia and Solern\'o proved that the sharper bound 
$ \ \deg g_i f_i \le n \, 2^{n+2} \ $ holds true \cite{SaSo95}. 
In fact more precise estimations are valid, we 
refer to the original papers for the exact statements. 

We note that such a degree bound allows us, given polynomials 
$f_1, \ldots, f_s \in k[x_1, \ldots, x_n]$, 
to determine whether the equation (1) is solvable or not, 
and in the case it is, to actually find a solution, as it reduces 
the original problem to the solving of a $k$--linear system 
of equations.

The study of this B\'ezout identity is the object of much research, 
due to both its 
theoretical and practical importance, mainly in the context of computational 
algebraic geometry and diophantine approximation. 
Thus it has been approached from many points of view and 
with different objectives. 
In this respect we refer to the 
research papers \cite{Amoroso96}, \cite{BeYg91}, 
\cite{Brownawell87}, \cite{CaGaHe89}, \cite{FiGa90}, \cite{GiHaHeMoMoPa96}, 
\cite{HaMoPaSo97}, \cite{KrPa96}, \cite{Philippon90}, \cite{SaSo95}, 
\cite{Shiffman89}, \cite{Smietanski94}.  
We also refer to the surveys 
\cite{BeSt91}, \cite{Pardo95}, \cite{Teissier91} for a broad introduction 
to the history of this problem, main results and open questions. 

\vspace{3mm} 

Let be given a Laurent polynomial 
$f = \sum_{i \in \ZZ^n } a_{i} \, x^{i}  \in k[x_1, \ldots, x_n, 
x_1^{-1}, \ldots, x_n^{-1}]$. The {\it support} 
of $f$ is defined as the set $\{ i  : a_{i} \neq 0\}$, that is, the 
set of exponents of the nonzero monomials of $f$. 
More generally, the support of a family of Laurent polynomials 
$f_1, \ldots, f_s $
is defined as the set of exponents of the nonzero monomials 
of any of the $f_i$  for $i=1,\ldots, s$. 

The {\it Newton 
polytope} $\N (f_1, \ldots , f_s)$ is defined as the convex hull of the
support of $f_1, \ldots, f_s$. 
Let $\rho $ denote the dimension of this polytope. 
Then the {\it unmixed volume }
${\cal U} (f_1, \ldots, f_s ) $ of the family 
of Laurent polynomials $f_1, \ldots, f_s$ is defined as $\rho !$ times 
the volume of the polytope $\N (f_1, \ldots, f_s)$. 

The degree of a polynomial is bounded by a nonnegative integer $d$ if 
and only if its Newton polytope is contained in $d$ times the standard 
simplex of $\RR^n$. 
Thus the notion of Newton polytope gives a sharper characterization of the 
monomial structure of a polynomial system than just degree. 
This concept is in the basis of sparse elimination theory. 
Within this theory, 
algorithms for elimination problems are designed trying to exploit the 
sparseness of  the involved polynomials, and sparseness is then usually 
measured in terms of the Newton polytope of  these polynomials. 
This is the point of  view introduced by Sturmfels in his foundational 
work  \cite{Sturmfels93} and followed, for instance, in \cite{CaEm96}, 
\cite{HuSt95}, \cite{Rojas97a}, \cite{Rojas97b}, 
\cite{VeVeCo96}. 

The sparse aspect in the Nullstellensatz has also been considered by Canny 
and Emiris, who
obtained a sparse effective Nullstellensatz but only for the case 
of $n+1$ {\it generic} $n$--variate Laurent polynomials 
\cite[Theorem 7.2]{CaEm96}. 

We obtain the following result, which in this context can be seen
as a bound 
for the sparseness of the output polynomials in 
terms of the sparseness of the input system. 

\begin{undef}{Theorem 1.} \label{thm1}
Let be given polynomials $ f_1, \ldots, f_s \in k[x_1, \ldots, x_n]$ 
without common zeros in $\AA^n$. 
Let ${\cal N}$ denote the Newton polytope of  the polynomials $x_1, \ldots, 
x_n, f_1, \ldots, f_s$,  and let $\U$ denote the 
unmixed volume of this polytope. 
Then there exist $g_1, \ldots, g_s  \in k[x_1, \ldots, x_n]$ satisfying
\[
1 = g_1 f_1+ \cdots + g_s f_s,
\]
with $ \ \N (g_i f_i) \subseteq \ n^{n+3} \, {\cal U} \ {\cal N}$ \ 
for $i=1, \ldots, s$. 
\end{undef}

Let us keep the notations of Theorem 1 and let 
$d$ denote the maximum degree of $f_i$ for $i=1, \ldots, s$. 
We readily derive from the previous result the degree bound
\[ 
\deg g_i f_i \le n^{n+3} \, d \ {\cal U}.
\]
We obtain from this the worst--case bound $\ \deg g_i f_i
\le n^{n+2} \, d^{n+1}$, as the unmixed volume of 
the polynomials $x_1, \ldots, x_n , f_1, \ldots, f_n$ is bounded by 
$ \ d^n$. We show however that our degree bound can 
considerably improve the usual one in the 
case when the input polynomial system is sparse 
and $d \ge n$ (Example \ref{exmp2.1}). 

We note that the naive notion of sparseness, based just on 
the counting of the number of nonzero monomials in each polynomial, does 
not yield better bounds for the degrees
than the usual ones, in view of the well known 
example of Mora, Lazard, Masser and Philippon \cite{Brownawell87}. 

We obtain an analogous result for the case of Laurent polynomials.

\begin{undef}{Theorem 2.} \label{thm2}
Let be given Laurent polynomials $ f_1, \ldots, f_s \in 
k[x_1, \ldots, x_n, x_1^{-1}, \ldots, x_n^{-1}]$ 
without common zeros in $(\overline{k}^*)^n$. 
Let $\N$ denote the Newton 
polytope of $f_1, \ldots, f_s$, and let 
$\U$ denote the unmixed volume of this polytope. 
Then there exist $a \in \ZZ^n$ and $g_1, \ldots, g_s  
\in k[x_1, \ldots, x_n, x_1^{-1}, \ldots, x_n^{-1}]$ 
satisfying 
\[
1 = g_1 f_1+ \cdots + g_s f_s,
\]
with $\ a \in \ n^{2n+3} \, {\cal U}^2 \ \N \ $  and $\ \N (g_if_i) 
\subseteq \ 
n^{2n+3} \, {\cal U}^2 \ {\cal N}-a \ $  for $i=1, \ldots, s$.
\end{undef}

\vspace{3mm}

The proof of both results is similar. 
It takes as its first step  the translation of the original 
system of equations over the affine space 
or the torus into a system of linear equations over an appropriate toric 
variety. 
The resulting system is then  solved by appealing to an effective 
Nullstellensatz for linear  forms  in a Cohen--Macaulay graded ring.  
This key lemma is proved following for the most part the lines 
of a previous paper \cite{Sombra96} which  in turn is based on previous 
work of Dub\'e  \cite{Dube93} and Almeida \cite{Almeida95}, although we 
introduce at this time some simplifications into the arguments. 
In particular we eliminate the use of estimations for the Hilbert function. 

As a by--product, we obtain an effective Nullstellensatz which holds 
not only  for 
linear forms, but for arbitrary homogeneous elements in a Cohen--Macaulay 
graded ring (Theorem \ref{thm1.1}). 

\vspace{3mm}

Besides we apply these arguments in two other situations. 
First we consider the usual effective Nullstellensatz and 
we obtain the following degree bound. 

\begin{undef}{Theorem 3.} \label{thm3}
Let be given polynomials $f_1, \ldots, f_s \in 
k[x_1, \ldots, x_n]$ without common zeros in $\AA^n$. 
Let $d_i$ denote the degree of $f_i$ for $i=1, \ldots, s$ and assume
that $d_1\ge \cdots \ge d_s$ holds. 
Then there exist $g_1, \ldots, g_s \in k[x_1, \ldots, x_n]$ satisfying
\[
1=g_1 f_1 + \cdots + g_s f_s,
\]
with $ \ \deg g_i f_i \le 2 \, d_s \, \prod_{j=1}^{\min\{ n, s \}-1}
d_j \ $ for $i= 1, \ldots , s$. 
\end{undef}

In the case when the polynomials $f_1, \ldots , f_s$ are quadratic
this estimation improves the bound 
$\ \deg g_i f_i \le n \, 2^{n+2} \ $ to the bound 
$\ \deg g_i f_i \le 2^{n+1}$, which is very close to the 
expected  $\ 2^n$. 

Finally we obtain another bound for the degrees in the 
Nullstellensatz. 
This bound depends
on the number of variables and the maximum degree of the input polynomials
and also on an additional parameter, called the 
{\it algebraic degree} of the input polynomial system, which we are going to 
introduce. 

Let $f_1, \ldots, f_s\in k[x_1, \ldots, x_n]$ be polynomials 
without common zeros in $\AA^n$. 
Let $ \lambda=(\lambda_{ij})_{ij}  
\in \overline{k}^{s\times s}$  be an arbitrary $s \times s$ matrix 
with entries in $\overline{k}$.  
We note by 
$h_i(\lambda)$ the linear combinations $\sum_j \, 
\lambda_{ij} \, f_j $ induced by the matrix $\lambda$ for $i=1,  \ldots, s$. 

Consider the set $\Gamma$ of $s\times s $ matrices 
such that for any $\lambda$ in $\Gamma$ the polynomials 
$h_1(\lambda), \ldots, h_{t-1}(\lambda)$ form a
regular sequence in $\overline{k}[x_1, \ldots, x_n]$ and 
$1\in (h_1(\lambda), \ldots, h_t(\lambda))$ for some 
$t=t(\lambda) \le \min\{n,s\}$.  
This set is nonempty, and indeed it contains a nonempty open set of 
$ \overline{k}^{s \times s}$.

For each $\lambda \in \Gamma$ and $i=1, \ldots, t-1$ we denote by  
$J_i(\lambda) \subseteq k[x_0, \ldots, x_n]$ the homogenization of 
the ideal $(h_1(\lambda ), \ldots, h_i(\lambda ))$. 
Then let $\delta(\lambda)$ denote the maximum degree 
of the homogeneous ideal $J_i(\lambda)$ for  $i=1 , \ldots, t-1$. 

The {\it algebraic degree} of the polynomial system $f_1, \ldots, f_s $ 
is defined as
\[
\delta(f_1, \ldots, f_s):= \min\, \{  \delta(\lambda):  \lambda \in \Gamma\}.
\]

This is the algebraic analogue of the notion of geometric degree of a system 
of equations 
of Giusti et al.\ \cite{GiHeMoMoPa95}, Krick, Sabia and Solern\'o 
\cite{KrSaSo96} and Sombra \cite{Sombra96}. 
We refer to Section 3 for a comparison between both notions.  

There have been obtained degree bounds for the polynomials in the 
Nullstellensatz which mainly depend on the geometric degree 
\cite{GiHaHeMoMoPa96}, \cite{KrSaSo96}, \cite{Sombra96}. 
We show that a similar bound
holds by replacing the geometric degree of the input polynomial
system by the algebraic degree. 

\begin{undef}{Theorem 4.} \label{thm4}
Let be given polynomials $f_1, \ldots, f_s \in k[x_1, \ldots, x_n]$ 
without common zeros in $\AA^n$. 
Let $d$ denote the maximum degree of $f_i$ for $i=1,  \ldots, s$ 
and let $\delta $ denote
the algebraic degree of this polynomial system. 
Then there exist $g_1, \ldots, g_s \in k[x_1, \ldots, x_n]$ satisfying 
\[
1=g_1 f_1 + \cdots + g_s f_s,
\]
with $ \ \deg g_i f_i \le \min\{n, s \}^2 \, d \, \delta \ $  
for  $i=1, \ldots, s$. 
\end{undef}

Let us denote by $d_i $  the degree of $f_i$ for $i=1, \ldots, 
s$ and assume that $d_1 \ge \cdots \ge d_s$ holds. 
Then the B\'ezout bound 
$ \ \delta(f_1, \ldots, f_s)\le d_s \, \prod_{i=1}^{\min\{ n,s\}-2} \, d_i \ $ 
holds true, and therefore we essentially recover from 
the previous result the known bounds for 
the degrees in the Nullstellensatz. 
The algebraic degree is bounded by the geometric degree, and so we also 
recover the known degree bounds in the Nullstellensatz which depend 
on the geometric degree. 
We show however that the algebraic degree is much smaller 
than the geometric degree in some particular instances, 
and by force, than the B\'ezout bound $\ d^{n-1}$ (Example \ref{exmp3.1}). 
We conclude that the obtained degree 
bound is much sharper in these cases than the known ones.

\vspace{3mm}

The outline of the paper is as follows. In Section 1 we obtain the effective 
Nullstellensatz for homogenous elements in a Cohen--Macaulay graded ring. In 
Section 2 we prove both Theorems 1 and 2 and we derive some of their 
consequences. In Section 3 we devote to the degree bounds in the usual 
Nullstellensatz.

\vspace{10mm}

% ----------- Seccion 0 -----------------------------------------------

\typeout{Section 0}

\noindent {\bf 0. Notations and Conventions}

\vspace{5mm}
\ 

We denote by  $k$ an infinite field with an algebraic closure denoted
by $\overline{k}$. 
All the rings to be considered in the sequel are Noetherian commutative, and
more precisely finitely generated $k$--algebras. 
The polynomial ring $k[x_0, \ldots, x_n]$ is alternatively denoted by $S$. 

As usual ${\PP}^n$ and $\AA^n$ denote the projective and the affine spaces 
of  $n$ dimensions over $\overline{k}$, respectively. 
A variety is not necessarily irreducible.

Let $J$ be a homogeneous ideal in the ring $S/I$. The {\it dimension}
of $J$ is defined as 
the  Krull dimension of the quotient ring $\ (S/I)/J \ $ and it is 
denoted by $\ \dim J$. 
The {\it degree} of $J$ is defined as $(\dim J-1)! $ times the 
principal coefficient of the Hilbert polynomial of the graded 
$k$--algebra $(S/I)/ J$. 

Let $\varphi : A\to B$ be a morphism between two rings, and let $I$ and $J$ 
be ideals of  $A$ and $B$ respectively. 
Then $I^e$ denotes
the extension of the  ideal $I$ to $B$ and $J^c$ denotes 
the contraction of the ideal $J$ to $A$. 
Given an element $\alpha \in A$,  we denote by $\overline{\alpha}$ the element 
$\varphi(\alpha)\in B $. 
We shall make use of this notation when it is clear from the context 
the morphism to which we refer.

\vspace{10mm}

% ----------- Seccion 1 -----------------------------------------------

\typeout{Section 1}

\setcounter{section}{1}
\setcounter{subsection}{0}

\noindent{\bf 1. An Effective Nullstellensatz over Cohen--Macaulay 
			 Graded Rings}

\vspace{5mm}

Let be given a Cohen--Macaulay homogeneous ideal
$I$ in the ring $k[x_0, \ldots, x_n]$. 
By this we mean that the quotient ring $k[x_0, \ldots, x_n]/I$ is 
a Cohen--Macaulay ring. 
Let $r$ denote the dimension of $I$. 
We also suppose that there is given a homogeneous element 
$p$ in $S/I$ which is not a zero--divisor. 
Let $\eta_1, \ldots, \eta_s \in S/I$ be homogeneous elements 
of degree one ---or for short, linear forms--- which define the empty variety 
in the open set $ \{ p \neq 0\} $ of $V(I)$.
In this situation, Hilbert's Nullstellensatz implies that $p$ belongs to 
the radical of the 
ideal $(\eta_1, \ldots, \eta_s)$, that is, 
$p \in \sqrt{(\eta_1, \ldots, \eta_s)} $. 
Equivalently we have  that 
$1$ lies in the ideal $(\overline{\eta}_1, \ldots, \overline{\eta}_s)$ spanned 
by $\overline{\eta}_1,   \ldots, \overline{\eta}_s$ in the ring $(S/I)_p$. 

We are going to give a bound for the minimal 
$D \in \NN$ 
such that $p^D$ falls into the ideal $(\eta_1, \ldots, \eta_s)$ 
(Main Lemma \ref{mainlem1.1}). This bound 
depends on the number of linear forms, and on the dimension and the 
degree of the ideal $I$. 
As a consequence of this result we derive an effective Nullstellensatz
for Cohen--Macaulay graded rings (Theorem \ref{thm1.1} and Corollary 
\ref{cor1.2}). 

\vspace{3mm}

Let $A$ be a ring and let  $\alpha_1, \ldots, \alpha_t$ be elements of $A$. 
Then $\alpha_1, \ldots, \alpha_t$ is called a {\it weak 
regular sequence} if  $\overline{\alpha}_i$ is not a zero--divisor
in the ring $A/ (\alpha_1, \ldots, \alpha_{i-1})$ for $i=1, \ldots, t$. 
We note that this definition differs from usual notion of regular 
sequence only in 
one point, namely that it allows $\overline{\alpha}_t$ to be a unit in 
$ A/ (\alpha_1, \ldots, \alpha_{t-1})$. 

\begin{lem}\label{lem1.1}
Let notations be as before. 
Then there exist linear forms 
$\zeta_1, \ldots, \zeta_t\in (\eta_1, \ldots, \eta_s) $ for some 
$t\le \min\{ r, s\}$ such that 
$\overline{\zeta}_1, \ldots, \overline{\zeta}_t $ is a weak regular
sequence in  $(S/I)_p$ and 
$1\in (\overline{\zeta}_1, \ldots, \overline{\zeta}_t)$ 
\end{lem}

\begin{proof}{Proof.}
Let $z_i:= \sum_{j=1}^s \, \lambda_{ij} \, \overline{\eta}_j$  be a
generic $k$--linear combination of 
$\overline{\eta}_1, \ldots, \overline{\eta}_s$  for $i=1, \ldots, s$. 
We obtain a maximal secant sequence
$z_1, \ldots, z_t$ for some $t\le s$, as the field $k$ is infinite. 
By this we mean that  
$ \dim (z_1, \ldots, z_i) = r-i-1$ for $ i=1, \ldots, t-1$ and that 
$1\in (z_1, \ldots, z_t)$. 
In particular $t\le \min \{ r,s\}$ holds true. 

The ring $(S/I)_p$ is Cohen--Macaulay as it is a localization of 
a Cohen--Macaulay ring. Then
$z_1, \ldots, z_t$ is a regular sequence in $(S/I)_p$. 

We take $\zeta_i := \sum_{j=1}^s \lambda_{ij} \, \eta_j$ for
$i=1, \ldots, t$. Thus $\zeta_i \in (\eta_1, \ldots, \eta_s)$ 
and $\overline{\zeta}_i=z_i$ for  $i=1, \ldots, s$ as required. 
\end{proof}

Thus we can suppose without loss of generality that $\overline{\eta}_1,  
\ldots, \overline{\eta}_s$ is  a weak regular 
sequence in $(S/I)_p$ and that $s\le r$ holds. 
We assume this from now on.                              
Next we are going to show that $\eta_1, \ldots, \eta_s$ can be replaced
by polynomials  of controlled degree which form a 
regular sequence in $S/I$ (Corollary \ref{cor1.1}). 

The following lemma is a generalization of \cite[Remark 4]{Heintz83}.

\begin{lem}\label{lem1.2}
Let be given a homogeneous unmixed ideal $K\subseteq k[x_0, \ldots, x_n]$ 
and points $\xi_1, \ldots, \xi_m\in \PP^n$ not lying in $V(K)$. Then there 
exists a homogeneous polynomial $g$ in $K$ 
such that $\deg g \le \deg K $ and that $ g(\xi_i)\neq 0 $ hold.
\end{lem}

\begin{proof}{Proof.}
The case when $K$ is a homogeneous prime ideal follows easily 
from \cite[Remark 4]{Heintz83}.

Let us consider the general case. For each associated prime ideal 
$P$ of $K$ we take a homogeneous polynomial $g_P$ such that 
$\ \deg g_P \le \deg P \ $ and $g_P(\xi_i ) \neq 0 $ for $i=1, \ldots, m$. 
Let $Q_P$ be the corresponding $P$--primary ideal in the decomposition of $K$.
Let $l(Q_P)$ denote the length of $Q_P$, that is, the length of 
$(S/Q_P)_P$ as a $S/P$--module. Then let 
\[
g:= \prod_P g_P^{\ l(Q_P)},
\]
where the product is taken over all the associated prime ideals of 
$K$. 
Then $g(\xi_i) \neq 0 $ for $i=1, \ldots, m$, and we have also 
that the polynomial $g$ lies in the ideal $K$ by 
\cite[Lemma 1]{BrMa80} and the degree bound $ \ \deg g \le 
\sum_P \, l(Q_P) \, \deg P = \deg K \ $ holds true by 
\cite[Proposition 1.49]{Vogel84}.
It follows that $g$ satisfies the stated conditions.
\end{proof}

In the sequel we shall denote by $J_i$ the contraction to the ring 
$S/I$ of the ideal $(\overline{\eta}_1, \ldots, \overline{\eta}_i) \subseteq 
(S/I)_p$ and by $\delta_i $ the degree of the 
homogeneous ideal $J_i$ for $i=1, \ldots, s$. 

\begin{cor} \label{cor1.1}
Let notations be as before. Then there exist homogeneous elements 
$h_1, \ldots, h_s\in S/I$ satisfying  the following conditions:

\begin{itemize}

\item[i)] $ h_i \equiv p^{c_i} \eta_i \ \ \ \mod J_{i-1}$
\ \ for some $c_i \ge 0$, 

\item[ii)] $h_1, \ldots, h_s$ is a regular sequence, 

\item[iii)] $\deg h_i \le \deg  J_{i-1} + \deg p-1$,

\end{itemize}

for $i=1, \ldots, s$.
\end{cor}

\begin{proof}{Proof.} 
We proceed by induction on $i$. 
By assumption $p$ is not a zero--divisor in $S/I$ so that the canonical 
morphism $S/I\to (S/I)_p$ is injective. 
The fact that $\overline{\eta}_1$ is not a zero--divisor in $(S/I)_p$ 
implies then that $\eta_1$ is not a zero--divisor in $S/I$. 

Now let $i\ge 2$ and assume that the elements $h_1, \ldots, h_{i-1}$ 
are already constructed. Let $H_{i-1}$ denote the ideal 
spanned by $h_1, \ldots, h_{i-1}$ in $ S/I$. Let  
$H_{i-1}=(\cap_j \, Q_j) \cap (\cap_l \, R_l)$ be the primary decomposition 
of $H_{i-1}$, with $p\notin \sqrt{Q_j}$ and $p\in \sqrt{R_l}$. 
Our objective is to find a homogeneous element  $h_i$ in $S/I $ 
not lying in any of the associated primary ideals
of $H_{i-1}$. 

We recall that the ideal $H_{i-1}$ has no imbedded component as
it is spanned by a regular sequence in a Cohen--Macaulay ring. 
On the other hand the ideal $J_{i-1}$ has the primary decomposition 
$ \cap_j \, Q_j$  and so it follows that 
$\ V(R_l) \subseteq \hspace{-4mm}/ \ V(J_{i-1}) \ $ holds for each $l$. 
We choose a point $\xi_l$ lying in $\ V(R_l)-V(J_{i-1}) \ $ and a 
homogeneous element $g\in J_{i-1}$ such that 
$\ \deg g \le \deg J_{i-1} \ $ and $g(\xi_l) \neq  0$  for each $l$.
The existence of $g$ is guaranteed by the previous lemma. 
By eventually multiplying 
$g$ with linear forms we can suppose without loss of generality 
that $\ \deg g = c_i \, \deg p +1 \ $ holds  
for some $ c_i \ge 0$. 
In particular we can assume that $ \ \deg g \le \deg J_{i-1} + \deg p -1 
\ $ holds. 
Finally we set 
\[
h_i:= a g + p^{c_i}  \eta_i 
\]
for some indeterminate scalar $ a \in k$. Then $h_i$ is homogeneous and
$h_i\equiv p^{c_i} \eta_i $  $\mod  J_{i-1} $ holds true. 
Therefore $h_i$ does not belong to 
$\sqrt{Q_j}$, as both $p$ and  $\eta_i$ are not zero--divisors 
modulo $J_{i-1}$. 
We have also that $ \ h_i(\xi_l)= a \, g(\xi_l) + (p^{c_i} \eta_i)(\xi_l) 
\neq 0 \ $ for a generic choice of $a$, which forces 
$h_i\notin \sqrt{R_l}$. 
\end{proof}

We fix the following notation. 
Let $h_1, \ldots, h_s\in S/I$ be the 
homogeneous polynomials introduced in  Corollary \ref{cor1.1},
and let $H_i:=(h_1, \ldots, h_i)$ and $L_i:=(\eta_1, \ldots, \eta_i) $ denote 
the homogeneous ideals successively generated by 
$h_1, \ldots, h_s$ and $\eta_1, \ldots, \eta_s$ respectively. 

Let us write $h_i=l_i + p^{c_i} \, \eta_i $ 
for some $l_i\in J_{i-1}$ and $c_i \ge 0$. 
Then set $\ \gamma_i:=\delta_{i-1}-\delta_i \ $, and let
$\ \lambda_i:= \sum_{j=1}^i (\gamma_j + c_j) \ $ and $\ \mu_i :=
\sum_{j=1}^i ( (i-j+1) \gamma_j + (i-j) c_j ) \ $ for $i=1, \ldots , s$. 

Given an ideal $K$ in $S/I$ we denote by  $K^u$
the unmixed part of $K$, that is, the unmixed ideal given as
the intersection of the primary components of $K$ of maximal 
dimension. 

\begin{lem} \label{lem1.4}
Let be given an element $q\in J_i$ for some $1\le i \le  s $. 
Then $ p^{\gamma_i} q \in (J_{i-1}, \eta_i)^u $.              
\end{lem}

\begin{proof}{Proof.}
Let $(\cap_j \, Q_j)  \cap  (\cap_l \, R_l)$ 
be the primary decomposition 
of the ideal  $(J_{i-1}, \eta_i)^u$, with $p \notin \sqrt{Q_j}$ and 
$p\in \sqrt{R_l}$. 
Then the ideal $J_i$ has the primary decomposition $\cap_j \, Q_j $. 
Let $K_i := \cap_l  \, R_l$ be the intersection of the
other primary components.
Then $K_i $ is an unmixed ideal which lies  in the hypersurface  $\{p=0\}$. 

The ideals $(J_{i-1}, \eta_i )^u $ and $(J_{i-1}, \eta_i)$  have the same
degree because they only differ in an ideal of codimension at least $i+1$.
The degree of $(J_{i-1}, \eta_i)$ is $\delta_{i-1}$, as $\eta_i$
is not a zero--divisor $\mod J_{i-1}$, and so the degree
of $K_i$ equals $\gamma_i= \delta_{i-1}-\delta_i$. 
Therefore $p^{\gamma_i}$ lies in the ideal $K_i$ 
\cite[Lemma 1]{BrMa80} and we 
conclude that $ \ p^{\gamma_i} q \in (\cap_j \, Q_j)  \cap  (\cap_l \, R_l)=
(J_{i-1}, \eta_i)^u \ $ as stated. 
\end{proof}

The following two statements (Lemmas \ref{lem1.5} and \ref{lem1.6}) 
are plain extensions of \cite[Lemmas 6.1 and 6.2]{Dube93}.

\begin{lem} \label{lem1.5}
Let be given an element $q\in J_i$ for some $1\le i\le s$. 
Then $p^{\lambda_i} q \in H_i$. 
\end{lem}

\begin{proof}{Proof.}
We proceed by induction on $i$. 
First  $p^{\gamma_1}q \in (\eta_1)^u$ by Lemma \ref{lem1.4}. 
We have also that $(\eta_1)^u= (\eta_1)$ and so the assertion 
is true for $i=1$. 

Let $i\ge 2$ and assume that the statement holds for $i-1$. 
By Lemma \ref{lem1.4}, $p^{\gamma_i} q \in (J_{i-1}, \eta_i)^u$,  that is, 
$p^{\gamma_i} \, q$ belongs to the intersection of the 
primary components of dimension $r-i$ of the ideal $(J_{i-1}, \eta_i)$. The 
intersection of the other primary components is an ideal
of codimension at least $i+1$.  
Then there exists a regular sequence $w_1, \ldots, w_{i+1}$ in this ideal, 
as $S/I$ is a Cohen--Macaulay ring. We have 
that $w_j \, p^{\gamma_i} \, q  \in (J_{i-1}, \eta_i)$ and so 
there exist $u_j\in J_{i-1}$ and $v_j \in S/I$  such that 
$ \ w_j\, p^{\gamma_i} \, q = u_j + v_j \eta_i \ $ for $j=1, \ldots, i+1$.
Then 
\[
w_j \, p^{\gamma_i +c_i} \, q = p^{c_i}\, u_j +  p^{c_i} \, v_j \, \eta_i=
p^{c_i}u_j + v_j (h_i - l_i) = (p^{c_i} u_j -v_j l_i ) + v_j h_i.
\]

Therefore $p^{\gamma_i+c_i}  u_j -v_j l_i \in J_{i-1}$ and
by the inductive hypothesis
$\ p^{\lambda_{i-1}}(p^{\gamma_i + c_i} u_j -v_j l_i) \ $ lies in the ideal 
$H_{i-1}$. 
Then $w_j \, p^{\lambda_i} \,q \in H_i$ holds for $j=1, \ldots, i+1$, 
as $\lambda_i = \lambda_{i-1} + \gamma_i -c_i$. 

The ideal $H_i$ is spanned by a regular sequence 
$h_1, \ldots, h_i $ and  so it is an unmixed 
ideal of dimension $r-i$. 
Thus for each associated prime ideal $P$ of $H_i$ 
there exists some $j$ such that $w_j\notin P$. 
We conclude that $p^{\lambda_i} q  \in H_i$.
\end{proof}

\begin{lem} \label{lem1.6}
Let be given an element $q \in J_i$  for some $1\le i \le s$. 
Then $p^{\mu_i} q \in L_i$. 
\end{lem}

\begin{proof}{Proof.}
We shall proceed by induction on $i$. 
The case $i=1$ follows in the same way as in the preceding lemma 
because $L_1 =H_1$ and $\mu_1 =\lambda_1$. 

Let $i\ge 2$. Then $p^{\lambda_i} q $ lies in $H_i$ by Lemma \ref{lem1.5}. 
Let us write 
$p^{\lambda_i} q  = u+v \, h_i$ for some  $u\in H_{i-1}$ and $v\in S/I$. 
Therefore $p^{\lambda_i} q - v \, h_i \in H_{i-1}$ and thus 
$p^{\lambda_i}q -  p^{c_i} \,  v \, \eta_i$ lies in the ideal $J_{i-1}$
because $H_{i-1}\subseteq J_{i-1}$ and 
$h_i \equiv p^{c_i} \, \eta_i \ \ \ \mod J_{i-1}$. This implies in turn that 
$p^{\lambda_i- c_i} q -  v \, \eta_i \in J_{i-1}$. 

From the inductive hypothesis we get that 
$p^{\mu_{i-1}}(p^{ \lambda_i- c_i}  \, q- v \, \eta_i $ 
lies in $ L_{i-1}$ and so $p^{\mu_{i-1} + \lambda_i- c_i}  q \in L_i$. 
The statement follows from the observation  
that $\mu_i= \mu_{i-1}+ \lambda_i-c_i$. 
\end{proof}

\begin{mainlem} \label{mainlem1.1} 
Let $I\subseteq k[x_0,\ldots, x_n]$ be a homogeneous 
Cohen--Macaulay ideal of dimension $r$. 
Let be given in addition 
a homogeneous element $p \in k[x_0, \ldots, x_n]/I $ 
which is not a zero--divisor and linear forms
$\eta_1, \ldots, \eta_s \in k[x_0,\ldots, x_n]/I$
such that $p$ lies in the radical of  the ideal $(\eta_1, \ldots, \eta_s)$. 
Then 
\[
p^D\in (\eta_1, \dots, \eta_s)
\]
holds, with $ \ D:= \min \{r,s\}^2   \deg I$. 
\end{mainlem}

\begin{proof}{Proof.}
By Lemma \ref{lem1.1} we can suppose without loss of generality that 
$\overline{\eta}_1. \ldots, \overline{\eta}_s$ is a weak regular sequence in 
$(S/I)_p$ and that $s\le r$.
After Lemma  \ref{lem1.6} it only remains to bound $\mu_s$. 
We make use of the estimations  $ \ \gamma_i , c_i \le \delta_{i-1} \ $ 
and we get the bound
\[
\begin {array}{rcl}
\mu_s&=& \sum_{j=1}^s ((s-j+1) \gamma_j +  (s-j) c_j) \\[2mm]
   &\le & \sum_{j=1}^s ((s-j+1) \delta_{j-1} +  (s-j) \delta_{j-1}) 
   \le s^2 \,  \deg I.
\end{array}
\]
\end{proof}

Compare with Caniglia, Galligo and Heintz 
\cite[Proposition 10]{CaGaHe89}, Smietanski \cite[Lemma 1.44]{Smietanski94},
and Sturmfels, Trung and Vogel 
\cite[Theorems 2.1 and 5.3]{StTrVo95}

\vspace{3mm}

The rest of the section is devoted to the extension of the previous 
result to the case when we are given homogeneous elements of arbitrary degree
instead of linear forms. 
First we  establish  some generalities about the Veronese imbedding. 

Let us denote by $N$ the integer $ \ \comb{n+d}{\ d} -1 \ $ and let 
$a_0, \ldots, a_N$ denote the exponents of the different monomials 
of degree $d$ in $S$. 
Let 
\[
v_d:\PP^n\to \PP^N, 
\msp 
x:= (x_0: \cdots: x_n)\mapsto (x^{a_0}: \cdots : x^{a_N})
\]
be the Veronese map. 
This is a regular morphism of projective varieties
and so its image is a closed subvariety of $\PP^N$. 
This variety is called
the Veronese variety and it is denoted by $v_{n,d}$. 
Let $I(v_{n,d})$ be its defining ideal and let us denote by $S(d):=
k[y_0, \ldots, y_N]/ I(v_{n,d})$ its homogeneous coordinate ring. 
The Veronese map induces then an inclusion of $k$--algebras 
$i_d: S(d) \hookrightarrow  S$ defined by 
$y_j\mapsto x^{a_j}$ for $j=0, \ldots, N$.

Let be given an ideal $J$ in $S$ and let $J(d)$  denote its contraction 
to the ring $S(d)$. 
Identifying the quotient ring $S(d)/J(d)$ with
its image in $S/J$ through the inclusion 
$i_d: S(d)/J(d) \hookrightarrow  S/J$ we obtain the decomposition  
in graded parts
\[
S(d)/J(d)= \oplus_j \, (S/J)_{d\, j}.
\]
Let $h_{J(d)}$ and $h_J$ denote the Hilbert functions of $J(d)$ and $J$
respectively. 
Then $h_{J(d)}(m)= h_J(d\, m)$ for  $m\in \NN$.
It follows that the ideals $J(d)$ and $J$ have the same dimension and  
that their degrees are related by the formula $ \ \deg J(d) =  
d^{\ \dim J-1 } \deg J $.

\begin{lem} \label{lem1.7}
Let $J$  be a homogeneous Cohen--Macaulay ideal in $S$ and let 
$J(d)$ denote its contraction to the ring $S(d)$. Then $J(d)$ is a 
Cohen--Macaulay ideal.
\end{lem}

\begin{proof}{Proof.}
Let us denote by $A$ and $B$ the quotient rings $S(d)/J(d)$ and $S/J$ 
respectively. 
We identify $A$ with its image in $B$ through the inclusion $i_d$. 

We are going to prove 
that $A$ is a Cohen--Macaulay ring. 
As $A$ is a graded ring
it suffices to exhibit a regular sequence of homogeneous elements of 
length equal to the dimension of $A$. 

Let $e$ denote the dimension of the ring $B$, which is also the dimension 
of $A$. 
Let $\beta_1, 
\ldots, \beta_e$ be a regular sequence in $B$ of homogeneous elements. 
Let $\alpha_i:= \beta_i^d$  for $i=1,\ldots, e$.
Then $\alpha_1, \ldots, \alpha_e$ are elements of $A$ which form a 
regular sequence in $B$, by \cite[Theorem. 16.1]{Matsumura86}. 
We affirm that they also form a maximal regular sequence in $A$. 
We only have to prove that $\alpha_i $ is not a zero--divisor 
in $A/(\alpha_1, \ldots,  \alpha_{i-1})$  for $i=1, \ldots, 
e$. Let $\zeta \in A$ be an element such that
$\zeta \, \alpha_i \in (\alpha_1, \ldots, \alpha_{i-1})$. 
Then there 
exist homogeneous elements $\zeta_1, \ldots, \zeta_{i-1} \in B$ such that
$\zeta= \zeta_1 \alpha_1 + \cdots + \zeta_{i-1}\alpha_{i-1}$ because 
$\alpha_1, \ldots, \alpha_{i-1}$ is a regular sequence in $B$. 
An easy
verification shows that $\zeta_1, \ldots, \zeta_{i-1}$ can  be chosen 
to lie in $A$, from where it follows 
that $\zeta \in (\alpha_1, \ldots, \alpha_{i-1})$. 
\end{proof}

\begin{thm} \label{thm1.1}
Let $I\subseteq k[x_0,\ldots, x_n]$ be a homogeneous Cohen--Macaulay ideal. 
Let be given in addition an homogeneous element 
$p \in k[x_1, \ldots, x_n]/I$ which is not a 
zero--divisor and homogeneous elements  
$f_1, \ldots, f_s \in k[x_0,\ldots, x_n]/I$
such that $p$ lies in the radical of the ideal 
$(f_1, \ldots, f_s)$. Let $r$ denote the dimension 
of $I$ and let $d$ the maximum degree of $f_i$ for $i=1, \ldots, s$.
Then 
\[
p^D\in (f_1, \ldots , f_s)
\]
holds, with $ \ D:= r^2 \, d^r \deg I $. 
\end{thm}

\begin{proof}{Proof.}
First we note that the zero locus in $V(I)$ of the polynomials
$\{ f_i\}_i $ equals the zero locus in $V(I)$ of the polynomials
$\{ x_j^{d- \deg f_i} f_i\}_{ij}$. 
We have also that $ x_j^{d- \deg f_i} f_i$ lies in the ideal 
$(f_1, \ldots , f_s)$ for all $i$ and $j$. 
Therefore we can suppose without loss of generality
that $ f_i$ is a homogeneous polynomial of degree $d$ for $i=1,  \ldots,  s$. 
We note however that the number of  input polynomials have been enlarged
in this preparative step. 

Let $i_d : S(d) \hookrightarrow S$ be the inclusion of $k$--algebras 
induced by the Veronese map and let $I(d)$ denote the contraction of the 
ideal $I$ to the ring $S(d)$. Then we have the inclusion  
$i_d : S(d)/I(d) \hookrightarrow S/I$ and the decomposition 
in graded parts $i_d(S(d)/I(d)) = \oplus_j (S/I)_{dj}$. 
We take a linear 
form $\eta_i\in S(d)/I(d)$ such that $i_d(\eta_i) = f_i$ for  $i=1, 
\ldots, s$, which exists as the inclusion $i_d$ is a bijection in degree one. 
We take also a homogeneous element $q \in S(d)/I(d)$ such 
that $i_d (q ) =p^d$. 

The map $v_d: V(I) \to V(I(d))$ is a dominant regular map
of projective varieties and so it is surjective.  
Therefore the zero locus of the linear forms $\eta_1, \ldots, \eta_s$ lies 
in the image of the zero locus of the polynomials $f_1, \ldots, f_s$. 
The common zeros of $f_1, \ldots, f_s$ 
lie in the hypersurface $\{p^d=0\}$ of $V(I)$ and we have in addition that 
$v_d( \{p^d=0\} )= \{ q =0\}$. 
Then the subvariety of $V(I(d))$ defined by $\eta_1, \ldots, \eta_s $ lies in the 
hypersurface $ \{ q =0\}$. 

By Lemma \ref{lem1.7} the ideal $I(d)$ is 
Cohen--Macaulay, and we have also that $q$ is not a zero--divisor 
modulo $I(d)$. Then we are in the hypothesis of 
the Main Lemma \ref{mainlem1.1}. As a consequence we obtain that
\[
q^{\, r^2 \, \deg I(d)}\in (\eta_1, \ldots, \eta_s)
\]
holds. 
Finally we apply the morphism $i_d$ to the previous expression
and we get that
\[
p^{d \, r^2  \, (d^{\, r-1} \deg I )} \in (f_1, \ldots , f_s)
\]
holds, as $\deg I(d)= d^{r-1} \deg I$. 
\end{proof}

\begin{cor} \label{cor1.2}
Let be given an ideal $I \subseteq k[x_1, \ldots, x_n]$. 
Assume furthermore that the homogenization $I^h$ of
the ideal $I$ in the ring $k[x_0, \ldots, x_n]$
is a Cohen-Macaulay ideal. Let 
$f_1, \ldots, f_s \in k[x_1, \ldots, x_n]$ be polynomials without 
common zeros in the
affine variety $V(I)$. Then there exist $g_1, \ldots, g_s \in k[x_1, 
\ldots, x_n]$ such that
\[
\msp 1\equiv  g_1 f_1 + \cdots +g_s f_s \msp (\mod I)
\]
holds,  with $\ \deg g_i  f_i \le (r+1)^2 \, d^{r+1} \, \deg I^h \ $ 
for $i=1, \ldots, s$. 
\end{cor}

\begin{proof}{Proof.}
By assumption the ideal $I^h$ is a Cohen--Macaulay homogeneous ideal of 
dimension $r+1$. We have also that $x_0$ is not a zero--divisor modulo 
$I^h$.  

Let $f_i^h$ denote the homogenization
of $f_i$ for $i=1, \ldots, s$. The homogeneous 
polynomials $f_1^h, \ldots, \f_s^h$ have no common zero 
in $V(\I)$  outside the hyperplane $\{x_0=0\}$.  
By Theorem \ref{thm1.1} there exist homogeneous  polynomials $v_1, \ldots, 
v_s \in S$  such that                            
\[
\msp x_0^{\, (r+1)^2 \, d^{r+1}} = v_1 f^h_1 + \ldots + v_s f^h_s \ssp
(\mod I^h)
\]
holds, with $ \ \deg v_i f^h_i  = (r+1)^2\, d^{r+1} \ $. The
corollary then follows by evaluating $x_0 : =1$. 
\end{proof}

Let notations  be as in  Corollary \ref{cor1.2}. In the case when $I$
is the zero ideal,  that is, 
in the setting of the classic effective Nullstellensatz, we get
the  degree bound 
\[
\deg g_i  f_i  \le (r+1)^2 \, d^{r+1}.
\]
									  
\vspace{10mm}

% ----------- Seccion 2 -----------------------------------------------

\typeout{Section 2}

\setcounter{section}{2}
\setcounter{subsection}{0}

\noindent{\bf 2.  The  Sparse Effective Nullstellensatz}

\vspace{5mm}

In this section we shall devote to the sparse effective Nullstellens\"atze
(Theorems 1 and 2) and to
the derivation of some of their consequences. 

\vspace{3mm}

First we introduce notation and  state some basic facts from polyhedral 
geometry and toric varieties. We refer to the 
books \cite{Fulton93} and \cite{Sturmfels96} 
for the proofs of these facts and for a more 
general background on these topics.

Let us be given a finite set of integer vectors $\A\subseteq \ZZ^n$.
The convex hull of $\A$ as a subset of $\RR^n$ is 
denoted by $\conv (\A)$. The cone over $\conv (\A)$ is denoted 
by $\pos (\A)$. 

The set $\A$ is {\it graded} if there exists an integer vector 
$\omega \in \ZZ^n $ 
such that $<a , \omega>=1$ holds for every $a\in \A$, that is, when the 
set $\A$ lies in an affine hyperplane which does not contain the origin. 

Let $\ZZ \A$ denote the $\ZZ$--module generated by $\A$. 
Let $\RR \A$ denote the linear space spanned 
by $\A$, so that $\ZZ \A $ is a lattice in $\RR \A$. Let $\rho $ denote the
dimension of this linear space. 
Then we consider the euclidian volume form 
in $ \RR A$, normalized in such a way that the lattice 
$\ZZ \A$ has covolume $\rho !$ or equivalently, such that each primitive 
lattice simplex has unit volume. The normalized {\it volume} 
$\Vol (\A)$  of the set $\A$ is defined as the volume 
of its convex hull with respect to this volume form. 

We get readily from the definition the bound
\[
\Vol (\A)  \le  \rho ! \, \vol (\conv (\A)), 
\]
where $\vol (\conv (\A))$ denotes the volume of the convex hull of 
$\A$ with respect to the usual
non--normalized volume form of $\RR^n$. 

Let $\NN \A$ denote the semigroup spanned by $\A$. This semigroup  
is always contained in the semigroup $ \ \pos (\A) \cap \ZZ \A $. 
The set $\A$ is said 
to be {\it normal} or {\it saturated} if the equality 
$ \ \NN \A = \pos (\A) \cap \ZZ \A \ $ holds. 

A polytope  $\P$ is said to  be {\it integral }
if it is the convex hull of a finite set of integer vectors. 

Let $\S$ be an integral simplex. Then $\S$ is said to be {\it unimodular}
if its interior contains no integral vector. Let $\{
s_1, \ldots, s_k\}$ be the set of vertices of $\S$. Then we have  
that $\S$ is unimodular if and only if the set of integer vectors
$\{s_2-s_1, \ldots , s_k-s_1\}$ is normal. 

Let $\P$ be an integral polytope.  A subdivision of 
$\P$ is said to be {\it unimodular} if it is formed by 
unimodular integral simplices. 

Given integral polytope $\P$ in $\RR^n$, we denote by $\A(\P)$ 
the set $\{1\}\times (\P \cap \ZZ^n)$, which is a graded set 
of integral vectors in $\ZZ^{n+1}$.
We note that the set $\A(\P)$ is normal in
the case when $\P$ admits a unimodular subdivision.
\vspace{3mm}

With respect to toric geometry, we shall follow the lines of 
\cite{Sturmfels96}. 
This point of view differs from the usual one in algebraic geometry. 
It is more combinatorial and suits better for our purposes. 

Let us be given again a finite set of integer vectors 
$\A=\{ a_1, \ldots, a_N\}$ in $\ZZ^n $. 
We associate to the set $\A$ the morphism 
\[
\varphi_\A : k[y_1, \ldots, y_N] \to 
k[x_1, \ldots, x_n, x_1^{-1}, \ldots, x_n^{-1}],
\ssp  
y_i\to x^{a_i}.                    
\]
The kernel of this map is a prime ideal $I_\A$  of $k[y_1, \ldots, y_N]$, 
called the {\it toric ideal} associated to the set $\A$. 
This ideal defines an {\it affine toric variety}  $X_\A$ as its 
zero locus in $ \AA^N $. This variety is irreducible and its 
dimension equals the rank of the $\ZZ$--module $\ZZ \A$. 

The   $k$--algebra $k[x_1, \ldots, x_n, x_1^{-1}, \ldots, x_n^{-1}]$ 
is the coordinate ring  of the torus  $(\overline{k}^*)^n$. Thus the map
$\varphi_\A$ induces a dominant map $(\overline{k}^*)^n\to X_\A$. The image 
of this map is called the {\it torus} $T_\A$ of the affine  toric variety
$X_\A$. This torus equals the open set $\{y_1 \cdots y_N  \neq 0\}$ of 
$X_\A$. 

The ideal $I_\A $ is homogeneous if 
and only if the set $\A$ is graded. In this case the set $\A$ 
defines a {\it projective toric variety } 
$Y_\A$ as the zero locus of the ideal $I_\A$ in the projective space 
$\PP^{N-1}$. The dimension of $Y_\A$ equals then the rank of 
$\ZZ \A$ minus one, and its degree equals the normalized volume of the set 
$\A$. 

Let $\A=\{ a_1, \ldots, a_N\}\subseteq \ZZ^n $ be a graded set. 
The intersection of the  projective variety 
$Y_\A$ with the affine chart 
$\{ y_i \neq 0\}\cong \AA^{N-1}$ equals the affine toric variety associated 
to the  set 
\[
\A -a_i:= \{a_1-a_i, \ldots, a_{i-1} -a_i, a_{i+1}- a_i , \ldots, a_N-a_i\}.
\]
In fact $Y_\A$ is irredundantly covered by the 
affine varieties $X_{\A-a_i}$, where $a_i$ runs over the vertices of the
polytope $\conv (\A)$. 

The $k$--algebra $k[y_1, \ldots, y_N]/ I_\A$ is canonically isomorphic to 
the semigroup algebra $k[\NN \A]$. 
This algebra is normal if and only if the set 
$\A$ is normal. 
We recall Hochster's theorem that 
the $k$--algebra $k[\NN \A]$ is a Cohen--Macaulay domain 
in the case when the set $\A$ is normal.

Let be given  an integral polytope $\P$ in $\RR^n$. 
This polytope determines a fan 
$\Delta_\P$ and an {\it abstract} complete toric variety 
$X_\P= X(\Delta_\P)$. 
This variety comes equipped with an ample Cartier divisor
$D_\P$. This Cartier divisor defines then a map 
$\varphi_\P:X_\P \to \PP^{N-1}$, where 
$N$ denotes the cardinality of the set $ \{ \P \cap \ZZ^n\}$. 
The image of this map is the projective 
variety $Y_{\A(\P)}$, where the set $\A(\P)$ is defined as before as 
$ \{1\} \times (\P \cap \ZZ^n)$ 
\cite[Section 3.4]{Fulton93}. The divisor $(n-1) D_\P$
is very ample \cite{EwWe91},  and so
the graded set $\A((n-1) \, \P)$ is normal.

\vspace{3mm} 

\begin{thm} \label{thm2.1}
Let be given polynomials $ p, f_1, \ldots, f_s \in k[x_1, \ldots, x_n]$ 
such that $p$ lies in the radical of the ideal $(f_1, \ldots, f_s)$. Let 
$\P$ be an integral polytope which contains the Newton polytope of the 
polynomials $ 1, x_1, \ldots, x_n, f_1, \ldots, f_s $. Assume furthermore 
that $\A(\P) $ is a normal set of integer vectors in $\ZZ^{n+1}$. 
Then there exist $D\in \NN$ and $g_1, \ldots, g_s  
\in k[x_1, \ldots, x_n]$ such that 
\[
p^D = g_1 f_1+ \cdots + g_s f_s
\]
holds, with $ \ D \le \  n! \,\min\{n+1, s\}^2 \, \vol (\P) \ $ and 
$ \ \N (g_i f_i) \subseteq \ (1+ \deg p) \, n!\, \min\{n+1, s\}^2 
\, \vol (\P)  \ \P \ $ for $i=1, \ldots, s$. 
\end{thm}

\begin{proof}{Proof.}
Let $ \B= \{ b_0, \ldots, b_N\}$ denote the set of integer vectors
$\P\cap \ZZ^n$, so that $\A(\P) = \{1\} \times \B$. Assume that 
$b_0 =(0, \ldots, 0)$. We consider the morphism of $k$--algebras
\[
\ssp
\psi: k[y_1, \ldots, y_N] \to k[x_1, \ldots, x_n],
\ssp  
y_i \mapsto x^{b_i}.
\]
The kernel of this morphism is the defining ideal $I_{\B-b_0}$
of the affine toric 
variety $X_{\B-b_0}$. This affine variety is the intersection
of the projective toric variety $Y_{\A(\P)}$ with the affine cart 
$\{ y_0 \neq 0\}$ of $\PP^ N$. 
In addition the map $\psi $ induces an isomorphism $\AA^n \to X_{\B-b_0}$.  

Let $\zeta_i$ be a polynomial of degree one in 
$k[y_1, \ldots , y_N]$ such that $\psi (\zeta_i)=f_i$ for 
$i=1, \ldots,  s$. 
We take also a polynomial $q$ in  
$k[y_1, \ldots, y_N]$ of degree less or equal 
to the degree of $p$ such that $\psi (q)=p$.  
Then $\zeta_1, \ldots, \zeta_s $ have no common zero in $X_{\B-b_0}$
outside the hypersurface $\{q = 0\}$. 

Let $\eta_1, \ldots, \eta_s, u$ denote the homogenization
of $\zeta_1, \ldots, \zeta_s, q $ in $ k[y_0, \ldots, y_N]$ respectively. 
Then the
linear forms $\eta_1, \ldots, \eta_s $ have no common zero in 
$Y_{\A(\P)}$ outside the hypersurface $\{y_0\, u =0\}$. 

By assumption the set $\A(\P)$ is normal, and so 
$I_{\A(\P)}$ is a Cohen--Macaulay prime homogeneous ideal of 
$k[y_0, \ldots, y_N] $ of dimension less or equal that $n+1$. 
We have also that $y_0 \, u $ is not a zero--divisor 
modulo $ I_{\A(\P)} $. 
Then we are in the hypothesis of the Main Lemma 
\ref{mainlem1.1}. 
Let $D$ denote the integer $ \ \min\{n+1, s\}^2 
\deg Y_{\A(\P)}$. 
We obtain that there exist homogeneous elements 
$\alpha_1, \ldots, \alpha_s \in k[y_0, \ldots, y_N]/I_{\A(\P)}$ 
of degree $ \ (1+\deg u) \, D -1 \ $ satisfying
\[
(y_0\, u)^D =\alpha_1 \eta_1 + \cdots + \alpha_s \eta_s.
\]
Finally we evaluate $y_0 := 1$ and we apply the map 
$\psi$ to the preceding identity. We get
\[
p^D =g_1 f_1+ \cdots +g_s f_s,
\]
where we have set $g_i(x): = \alpha_i(1, x^{b_1 } , \ldots, x^{b_N})$ 
for $i=1, \ldots, s$. 
We have the estimations $\ \deg u \le  \deg p \ $ and 
$\ \deg Y_{\A(\P)} \le  n! \, \vol (\P) $.
We conclude that $\ D\le \, n! \, \min\{n+1, s\} \, \vol (\P) \ $ and 
that the polytope $\N(f_i g_i)$ is contained in 
$\ ( (1+\deg p) \, n!\, \min\{n+1, s\}^2 \, \vol (\P) ) \ \P \ $ for 
$i=1, \ldots, s$. 
\end{proof}

We derive  from the previous theorem the following degree  
bound.                              

\begin{cor} \label{cor2.1}
Let  notations be  as  in Theorem \ref{thm2.1}. 
Let $d$ denote the maximum degree of the  polynomials $f_i$ for  
$i=1,  \ldots,  s$. 
Then there exist $D\in \NN$ and $g_1, \ldots, g_s  
\in k[x_1, \ldots, x_n]$ such that 
\[
p^D = g_1 f_1+ \cdots + g_s f_s
\]
holds, with $\ D \le \  n! \,\min\{n+1, s\}^2 \, \vol (\P) \ $ and 
$\ \deg g_i f_i \le \ d\, (1+ \deg p) \, n!\, \min\{n+1, s,\}^2 
\, \vol (\P) \ $ for $i=1, \ldots, s$. 
\end{cor}

\begin{proof}{}
\end{proof}

We are going to show with an example that this degree bound can be much more 
precise than the usual one in the case when the
input system is sparse. 
First we need the following auxiliary result.

\begin{lem} \label{lem2.1}
Let $\P_d$ denote the integral polytope $\{ x=(x_1, \ldots, x_n)\in \RR^n: 
x_1, \ldots, x_n \ge 0, \  
x_1 + \cdots + x_{n-1} \le 1,  \ x_n \le d\}$ for some $d\in \NN$.
Then $\P_d$ admits a unimodular subdivision. 
\end{lem}

\begin{proof}{Proof.}
Let us denote by $e_1, \ldots , e_{n-1}$ 
the standard unit vectors of $\RR^{n-1}$. 
Let $\S$ denote the simplex of $\RR^{n-1}$ with vertices 
$(0, \ldots, 0),  e_1, \ldots, e_{n-1}$, so that 
$\P_d=\S \times [0, d]$. 
We have the  
subdivision $\ \P_d = \cup_{i=0}^{d-1} (\P_1 + i \, e) \ $, where 
$e\in \RR^n$ denotes the vector $(0, \ldots, 0,1)$.
Thus it suffices to show that $\P_1$ admits a unimodular subdivision. 

Let $\alpha_i$, $\beta_i \in \RR^n$ denote the vectors  $e_i\times \{0\}$,
$e_i \times\{1\}$ respectively, for $i=1, \ldots, n-1$,  and 
$(0, \ldots, 0, 0) $, $(0, \ldots, 0, 1)$ for $i=n$. 
Let $\S_i$ be the integral simplex determined by the vectors 
$\alpha_1, \ldots, \alpha_i, 
\beta_i, \ldots, \beta_n$ for $i=1, \ldots, n$. 

It is easily checked that the simplices $\S_i$ and $\S_j$
intersect in a proper face for $i\neq j$, so that 
$\S_1, \ldots, \S_n$ are essentially 
disjoint. We have also that $\vol (\S) =1/(n-1)!$ and that 
$\vol (\S_i) =1/n!$ for $i=1, \ldots, n$, 
and so $\S_1, \ldots, \S_n$  form a unimodular subdivision of $\P_1$
as desired.
\end{proof}

\begin{exmp} \label{exmp2.1}
{\rm 
Let $f_1, \ldots, f_s \in k[x_1, \ldots, x_n]$ be polynomials 
without  common  zeros in  $\AA^n$. 
Assume in addition that 
$ \ \deg_{x_1, \ldots, x_{n-1}} f_i \le 1 \ $ and that 
$\ \deg_{x_n} f_i \le d \ $ holds for 
$i=1, \ldots, s$. 
Let $\P_d $ denote as before the polytope 
$\{ x\in  \RR^n: x_1,  \ldots,  x_n \ge 0, \, 
x_1 + \cdots + x_{n-1} \le 1, \, x_n \le d\}$. 
Then $\P_d$ contains the 
Newton  polytope  of  the polynomials 
$1, x_1, \ldots, x_n, f_1, \ldots, f_s$ and the 
set $\A(\P)$ is normal by the preceding lemma.
Thus we are in the hypothesis
of Corollary \ref{cor2.1} and we can conclude that
there exist $g_1, \ldots, g_s \in k[x_1, \ldots, x_n]$ such that
\[
1=g_1 f_1 + \cdots + g_s f_s
\]
holds, with $ \ \N (g_i f_i) \subseteq d \, n \, 
\min\{ n+1, s \}^2 \, \P_d \ $ 
for $i=1, \ldots,  s$, as the volume of $\P_d$ equals $d/(n-1)!$.
In particular we get the degree bound 
$\ \deg g_i  f_i \le (n+1)^3\, (d+1)^2$, which is much sharper than 
the estimation 
$\ \deg g_i f_i \le (d+1)^n \ $ which follows from direct application 
of the usual degree bound.
}
\end{exmp}

Let notations be again as in Theorem \ref{thm2.1}. 
Let $\N$ denote the Newton polytope 
of the polynomials $1, x_1, \ldots,  x_n,
f_1, \ldots, f_s$ and let $\U$ denote the unmixed volume of this polytope. 
Assume that $n\ge 2$. 
In this situation we can then take the  polytope 
$\P$ to be  $(n-1) \,  \N$. 
Then we get the bounds 
\[
D\le n^{n+2} \, \U, 
\msp 
\N (g_i f_i)  \subseteq ((1+\deg p)  \, n^{n+3} \, \U) \ \N.
\]
It is easy to check that these bounds hold also when $n=1$. 
Thus Theorem 1 follows from this observation in the particular case 
$p=1$. We observe that in this case the condition
$0 \in \P$ is redundant. 

\vspace{3mm}

We obtain a similar result in the case of Laurent polynomials. 

\begin{thm} \label{thm2.2}
Let be given Laurent polynomials 
$ p, f_1, \ldots, f_s \in k[x_1^{-1}, \ldots, x_n^{-1}, x_1, \ldots, x_n]$ 
such that $p$ lies in the radical of the ideal 
$(f_1, \ldots, f_s)$. Let $\P$ be an integral polytope which contains the 
Newton polytope  of $ p, f_1, \ldots, f_s$. Let $\rho$ denote  its  
dimension. Assume furthermore that $\A(\P) $ is 
a normal set of integer vectors in $\ZZ^{n+1}$.  
Then there exist $D \in \NN$, $a\in \ZZ^n$ and $g_1, \ldots, g_s  
\in k[x_1, \ldots, x_n, x_1^{-1}, \ldots, x_n^{-1}]$ such that 
\[
p^D = g_1 f_1+ \cdots + g_s f_s
\]
holds, with $\ D \le \ \rho! \,\min\{n+1,s\}^2 \, \vol (\P) \ $, 
$\ a \in (\rho! \,\min\{n+1, s\} \, \vol (\P))^2 \P \ $
and $\ \N(g_i f_i) \subseteq \ 
( \rho! \, \min\{n+1, s\} \, \vol (\P) )^2 \ \P - a \ $ for $i=1, \ldots, s$. 
\end{thm}

\begin{proof}{Proof.}
As before, we denote by $ \B= \{ b_0, \ldots, b_N\}$ the set 
of integer vectors $\P\cap \ZZ^n$.
Assume for the moment that $b_0 =(0, \ldots, 0)$. 
We consider the morphism 
\[
\psi:  k[y_1, \ldots, y_N] \to k[x_1, \ldots, x_n, x_1^{-1}, 
		 \ldots, x_n^{-1}], 
\ssp  
y_i \mapsto x^{b_i}.         
\]
The kernel of this morphism is the defining ideal $I_{\B-b_0}$
of the affine toric variety $X_{\B-b_0}$. 
Let $T$ denote  the  torus of this  toric  variety. Then we have that 
$X_{\B-b_0}$ equals the intersection of the  projective variety
$Y_{\A(\P)}$ with the  affine  cart  $\{y_0 \neq  0\}$ of  $\PP^N$,
and that  $T$ is also the torus  of $Y_{A(\P)}$. We recall that this  torus 
equals the open set 
$\{y_0  \cdots y_N  =0\}$ of  $Y_{A(\P)}$.  

The  map $\psi$ induces  a  surjection $(\overline{k}^*)^n  \to T$. 
Let $\zeta_1, \ldots, \zeta_s, q$ be elements of degree one in 
$k[y_1, \ldots , y_N]$ such that $\psi (\zeta_i)=f_i$
for $i=1, \ldots, s$ and $\psi (q) = p$. 
Then $\zeta_1, \ldots, \zeta_s $ have no common zero in  $T$ 
outside the hyperplane $\{q  =0\}$. 

Let $\eta_1,  \ldots, \eta_s, u$ denote the homogenization of 
$\zeta_1, \ldots, \zeta_s, q$ in  
$k[y_0, \ldots, y_N]$ respectively.
Then the linear forms $\eta_1, \ldots, \eta_s$ have no
common zero in  $Y_{\A(\P)}$ outside the 
hypersurface $\{  y_0 \cdots y_N \, u=0\}$. 

Let $V(\eta_1, \ldots, \eta_s)$ denote the subvariety of $Y_{\A(\P)} $ 
defined by the linear forms $\eta_1, \ldots, \eta_s $.  
By B\'ezout inequality, 
the number of irreducible components of $V(\eta_1, \ldots, \eta_s)$ 
does not exceed the degree of $Y_{\A (\P)}$. 
Let  us denote  by $\delta $ the degree of $Y_{\A(\P)}$,  
so that $ \ \delta \le  \rho! \, \vol (\P) \ $ holds. 
In our situation this implies that 
$V(\eta_1, \ldots, \eta_s)$ lies in the union of at most $\delta$ 
hyperplanes. 
These hyperplanes are defined by variables $y_{i_1}, \ldots,  y_{i_l}$, 
and eventually also by the linear form $u$, depending on whether $\eta_1,  
\ldots, \eta_s$ have a common  zero  in $T$ in the hyperplane $\{u =0\}$ 
or not. Let $\Pi$ denote the product of these equations, which is a 
polynomial of degree less or equal that $\delta$. 

By assumption the  set  $\A(\P)$ is normal and  so 
$I_{\A(\P)}$ is a Cohen--Macaulay prime homogeneous ideal of 
$ k[y_0, \ldots, y_N]$ . We have also that $\Pi$ is not  a  
zero--divisor modulo this ideal. 
Thus we are again in  the hypothesis of the Main Lemma 
\ref{mainlem1.1}. Let $E$ denote the integer 
$\ \min\{n+1, s\}^2 \, \deg Y_{\A(\P)}$. 
Then there exist homogeneous elements 
$\alpha_1, \ldots, \alpha_s \in k[y_0, \ldots, y_N]/I_{\A(\P)}$ 
of degree $ \ \deg \Pi \ E -1 \ $ such that
\[
\Pi^{E}= \alpha_1 \eta_1 + \cdots + \alpha_s \eta_s 
\]
holds. We evaluate $y_0 :=1$ and we apply the map $\psi $ to the preceding  
identity. We get
\[
p^D=g_1 f_1+ \cdots +g_s f_s,
\]
where we have set $g_i(x):= (x^{b_{i_1}} \cdots x^{b_{i_l}})^{-1}  \,  
\alpha_i (1, x^{b_1},  \ldots, x^{b_N})$ for $i=1, \ldots, s$  and
$D:=E$  in the case when $u$ appears as a factor of $\Pi$ and 
$D:=1$ in the other case. 
Then  $\ D \le \rho! \, \min\{n+1, s\}^2 \, \vol(\P) \ $ holds and 
the polytope $\N(g_if_i) $ is contained in 
$\ (\rho!\, \vol(\P) \, E -1 )\ \P - (b_{i_1} + \cdots + b_{i_l}) \ $
for $i=1, \ldots,  s$. We have that $ \ \deg \Pi  \le  \deg Y_{\A(\P)}\le 
\rho! \, \vol (\P) \ $ and that $i_1  +  \ldots  + i_k  \in  
\deg Y_{\A(\P)} \ \P$. 

Now we  consider the general case. Let  $b_0$ be any integer vector in $\P$,  
and let $\Q$ denote the polytope $\P-b_0$. 
By  the  previous considerations there exist $D\in \NN$, 
$a_0\in \ZZ^n$ and $g_1, \ldots, g_s \in k[x_1, \ldots, x_n, 
x_1^{-1} , \ldots, x_n^{-1}]$ such that
\[
p^D=g_1 f_1+ \cdots +g_s f_s
\]
holds, with $\ D \le \ \rho! \, \min\{n+1, s\}^2 \, \vol (\Q) \ $,   
$\ a_0\in  \rho! \, \vol (\Q) \ \Q \ $   and  
$\ \N(g_if_i) \subseteq \ (\rho! \, \min\{n+1, s\} \, 
\vol (\Q))^2\ \Q- a_0 \ $
for $i=1, \ldots, s$. 

Let  $a$ be the integer vector $ \ a_0+(\rho! \, \min\{n+1, s\} \, \vol (\P))^2 
\, b_0 $. Then $a$ lies in the polytope  
$\ (\rho! \, \min\{n+1, s\} \, \vol (\P))^2 \ \P \ $ and we have also that 
$\ \N(g_i f_i) \subseteq \ 
( \rho! \, \min\{n+1, s\} \, \vol (\P) )^2 \ \P - a \ $ holds 
for $i=1, \ldots, s$
as stated. 
\end{proof}

Let notations be as in Theorem \ref{thm2.2}. 
Let $\N$ denote the Newton polytope 
of $p, f_1, \ldots, f_s$ and let $\U$ denote the unmixed volume of this 
polytope. Assume in addition that $n\ge 2$. 
In this situation we can then take the  polytope 
$\P$ to be  $(n-1) \,  \N$. 
We get the bounds 
\[
D\le n^{n+2} \, \U, 
\msp 
\N (g_i f_i)  \subseteq (n^{2n+3} \, \U) \ \N -a.
\]
for some $\ a \in  (n^ {2n+3} \, \U) \, \N$. 
As before, it is  easy to verify that 
the same bounds hold also when $n=1$. 
Thus Theorem  2 follows from 
this observation in the particular  case $h=1$. 

Let be given a rational function $q\in k(x_1, \ldots, x_n)$ and let 
$q=f/g$ be a representation of $q$ as the quotient of two polynomials
without common factors. 
Then the {\it degree} of $q$ is defined as 
$ \ \deg q:= \max \{\deg f, \deg g\}$. 

We derive  from Theorem \ref{thm2.2} the following degree  
bound.                              

\begin{cor} \label{cor2.2}
Let  notations be  as  in Theorem \ref{thm2.2}. 
Let $d$ denote the  maximum  degree of the Laurent polynomials $f_i$ for  
$i=1,  \ldots,  s$. 
Then there exist $D\in \NN$ and $g_1, \ldots, g_s  
\in k[x_1, \ldots, x_n, x_1^{-1}, \ldots, x_n^{-1}]$ such that 
\[
p^D = g_1 f_1+ \cdots + g_s f_s
\]
holds, with $ \ D \le \ \rho! \,\min\{n+1,s\}^2 \, \vol (\P) \ $, 
$ \ a \in (\rho! \,\min\{n+1, s\} \, \vol (\P))^2 \P \ $
and $ \ \deg (g_i f_i) \le \ 
d\, ( (1+\deg p ) \rho! \, \min\{n+1, s\} \, \vol (\P) )^2 \ $ 
for $i=1, \ldots, s$. 
\end{cor}

\begin{proof}{}
\end{proof}

\vspace{10mm}

% ----------- Seccion 3 -----------------------------------------------

\typeout{Section 3}

\setcounter{section}{3}
\setcounter{subsection}{0}

\noindent{\bf 3. Improved Bounds for the Degrees in the 
			   Nullstellensatz}

\vspace{5mm}

In this section we consider the degree bounds in the Nullstellensatz.
We shall apply the methods used in Section 1 in a direct way ---without 
any reference to  the Veronese map---                       
in the setting of the classic effective Nullstellensatz. The proof 
follows closely the same lines and so we shall skip some verifications
in order to avoid unnecessary repetitions. 

\vspace{3mm}

Assume that we are given homogeneous polynomials 
$f_1, \ldots, f_s $ in  $k[x_0, \ldots, x_n]$ without common zeros at finite 
distance. In this situation we are going to give a bound for the 
minimal $D\in \NN$ such that 
$x_0^D \in (f_1, \ldots, f_s)$. 

We shall assume without loss of generality that 
$s\le n+1$ and  that $\overline{f}_1, \ldots, 
\overline{f}_s $  is a weak regular sequence in $k[x_0, \ldots, x_n]_{x_0}$.
Let $d_i$ denote the degree of $f_i$ for $i=1,  \ldots, s$. 
We shall also suppose that $d_2 \ge \cdots 
\ge d_s$ and that $d_s \ge d_1$ hold. 
As before these polynomials can 
be obtained 
as linear combinations of the original polynomials, 
eventually multiplied by powers of $x_0$.

Let us denote by  $J_i$ the contraction to the ring $S$ of the ideal
$(\overline{f}_1, \ldots, \overline{f}_i) \subseteq S_{x_0}$ 
 for $i=1, \ldots, s$. We make the convention $J_0:=(0)$.

\begin{lem} \label{lem3.1}
Let  notations be as before. Then 
there exist homogeneous polynomials $h_1, \ldots, h_s\in 
k[x_0, \ldots, x_n]$ satisfying  the following conditions:

\begin{itemize}

\item[i)] $h_i\equiv x_0^{c_i}f_i \ \ \ \mod J_{i-1}$
 for some $c_i \in \NN$, 

\item[ii)] $h_1, \ldots, h_s $ is a regular sequence, 

\item[iii)]$\deg h_i \le \max \, \{\deg J_{i-1}, \deg f_i \}$,

\end{itemize}
for $i=1, \ldots, s$.
\end{lem}

\begin{proof}{}
\end{proof}

We introduce the following notation. 
Let $\delta_i$ denote the degree of the  
homogeneous ideal $J_i$ for $i=0, \ldots, s$. 
We recall the B\'ezout bound $ \ \delta_i \le \prod_{j=1}^i \, d_j$. 
Then we denote by $\gamma_i$ 
the integer $ \ d_i \, \delta_{i-1}-\delta_i$ for $i=1, \ldots, \min\{n, s\} \ $
and $ \ \gamma_{n+1}:=\delta_n+ d_{n+1}-1$. 
We also let 
$\ \delta:=\max \{\delta_i : i=1, \ldots, s-1\} \ $ and 
$\ d:=\max \{d_i : i=1, \ldots, s-1\}$. 
 
We recall that given an ideal $I$ in $S$ we denote by  $I^u$ its  unmixed 
part.

\begin{lem} \label{lem3.2}
Let be given a polynomial $q\in J_i$ for some $1\le i\le s$. 
Then  $x_0^{\gamma_i} \, q \in (J_{i-1}, \eta_i)^u$. 
\end{lem}

\begin{proof}{Proof.}
The case $i\le n$ is exactly as in Lemma \ref{lem1.4}. Thus we only
consider the case $i=n+1$. 

The ideal $J_n$ is has dimension one and its degree is 
$\delta_n$. 
Then $(J_n, f_{n+1})_m=S_m$ for $m\ge \delta_n + d_{n+1}-1 $ 
as $f_{n+1}$ is not a zero--divisor 
modulo $J_n$ \cite[Theorem. 2.23]{Sombra96}. 
It follows that $x_0^{\gamma_{n+1}}\in (J_n, f_{n+1})$ and in particular 
$x_0^{\gamma_{n+1}}\, q \in (J_n, f_{n+1})^u$. 
\end{proof}

Now let $h_1, \ldots, h_s$ be the homogeneous polynomials introduced in 
Lemma  \ref{lem3.1}. We set  
$\ \mu_i:=\sum_{j=2}^i ((i-j+1) \, \gamma_j +  (i-j) \, c_j) \ $ for $i=1, 
\ldots, \min\{n, s\}$ 
and $\ \mu_{n+1}:= \mu_n + \gamma_{n+1}$, where $c_i$ denotes the integer  
$\deg h_i- \deg f_i$. 

We denote by $L_i$ the homogeneous ideal $(f_1, \ldots, f_i) $ 
for $i=1, \ldots, s$. 

\begin{lem} \label{lem3.3}
Let be given a polynomial $q \in J_i$ for some $1\le i\le s$. Then 
 $x_0^{\mu_i} q \in L_i$.
\end{lem}

\begin{proof}{Proof.}         
The case $i\le n$ is exactly as in Lemma \ref{lem1.6}. Thus we only consider the 
case $i=n+1$. 

By the previous lemma $x_0^{\gamma_{n+1}} \, q \in (J_n, f_{n+1})^u
=(J_n, f_{n+1})$ and so $x_0^{\gamma_{n+1}} \, q - u\, f_{n+1} 
\in J_n$ for some polynomial $u\in S$. 
We apply then the inductive hypothesis
and we obtain that $x_0^{\mu_n}\, (x_0^{\gamma_{n+1}} \, q - u\, f_{n+1} )
\in L_n$ from where it follows that $ x_0^{\mu_{n+1}}\, q \in L_{n+1}$. 
\end{proof}

Thus it only remains to bound $\mu_s$. We shall be concerned with
two different types of bounds. One depends as usual
on the number of variables and on the degrees of the input polynomials, 
and the other depends also on the degree of some ideals 
associated to these polynomials. 

\begin{lem} \label{lem3.4}
Let notations be as before. Then 
$ \ \mu_s \le \min \{ n,s\}^2  \, d \,  \delta $. In the case 
when $\ \deg f_i \ge 2 \ $ 
for $i=1, \ldots, s$ we have that  $ \ \mu_s \le 2 \, 
\prod_{j=1}^{\min \{ n, s\}} \, d_j$. 
\end{lem}

\begin{proof}{Proof.}
We decompose the integer $\mu_s$ in two terms and we estimate them 
separately. First we consider the term 
$\sum_{j=2}^s \, (s-j) \, c_j $. We have that  $c_i \le 
\max\{\delta_{i-1}- d_i, 0\}$. 
In particular $c_2 =0$ as $\delta_1 =d_1$ and 
$d_1 \le d_2$. 
Then 
\[
\begin{array}{rcl}
\sum_{j=2}^s \, (s-j) \, c_j &\le&  
\sum_{j=3}^{s-1} \, (s-j) \, (d_1 \cdots d_{j-1} -d_j) \\[2mm]
&\le & (\sum_{j=3}^{s-1} (s-j)/d_j \cdots d_{s-2})\, d_1 \cdots d_{s-2} - 
\sum_{j=2}^{s-1} \, (s-j) \, d_j \\[2mm]
&\le& 4 \, d_1 \cdots d_{s-2} - d_{s-1},
\end{array}
\]
under the assumption $d_i \ge 2$ for $i=1, \ldots, s$. 
We have also $\sum_{j=2}^{s-1} \, (s-j) \, c_j \le  
\sum_{j=2}^{s-1} \, (s-j) \  \delta 
= \frac{1}{2}\, (s-2)(s-1) \, \delta$. 

Now we estimate the other term. We consider first the case $s\le n$. Then 
\[
\begin{array}{rcl}
\sum_{j=2}^s \, (s-j+1) \, \gamma_j &=& 
\sum_{j=2}^s \, (s-j+1) \, (d_j \, \delta_{j-1} -\delta_j) \\[2mm]
&=& (s-1) \, d_2 \, \delta_1 +
\sum_{j=3}^s \, ((s-j+1) \, d_j - (s-j)) \, \delta_{j-1} -\delta_n\\[2mm]
&\le& d_1 \cdots d_s - \delta_s,
\end{array}
\]
from where we obtain the bound $\mu_s = \sum_{j=2}^s \, (s-j+1) \, \gamma_j 
+\sum_{j=2}^{s-1} \, (s-j) \, c_j 
\le (d_1 \cdots d_s - \delta_s) + (4 \,  d_1 \cdots d_{s-2} - d_{s-1}) 
\le 2 \, d_1 \cdots d_s$. 

In the case $s=n+1$ we have that $\ \mu_{n+1}=\mu_n+ 
\gamma_{n+1} \ $ from where it follows that $ \ \mu_{n+1} \le 
(2\, d_1 \cdots d_n - \delta_n - d_{n-1}) + (\delta_n + d_{n+1}-1)
\le 2\, d_1 \cdots d_n$.

On the other hand we have also the estimation  
$\sum_{j=2}^s \, (s-j+1) \, \gamma_j \le \frac{1}{2}\, (s-1)\, s 
\, d \, \delta $ from where we conclude that $ \ \mu_s \le  
\frac{1}{2}\, (s-1)\, s \, d \, \delta + \frac{1}{2} \, (s-2)\, (s-1) 
\, \delta \le (s-1)^2 \, d \, \delta \ $ holds, as stated. 
\end{proof}

\begin{thm} \label{thm3.1}
Let be given homogeneous polynomials $f_1, \ldots, f_s \in 
k[x_0,\ldots, x_n]$ such that $x_0 $ lies in the radical of the ideal 
$(f_1, \ldots, f_s)$. 
Let $d_i$ denote the degree of $f_i$ for $i=1, \ldots, s$ and assume
that $d_1\ge \cdots \ge d_s$ holds. 
Then 
\[
x_0^D\in (f_1, \ldots, f_s)                       
\]
holds, with $ \ D:=2 \, d_s\, \prod_{i=1}^{\min\{ n,s\}-1 }\, d_i$. 
\end{thm}

\begin{proof}{Proof.}
After Lemmas \ref{lem3.3} and \ref{lem3.4} it only remains to consider the 
case when some $f_i$ has degree one. 

By assumption $f_1, \ldots, f_s$ are ordered in such a way 
that $d_1 \ge \cdots \ge d_s$ holds. 
Let $r$ be maximum
such that $d_r\ge 2$, so that the polynomials 
$f_{r+1}, \ldots, f_s$ have all degree one. 
We can assume without loss
of generality that they are $k$--linearly independent. 
We can also suppose that neither $1$ nor $x_0$ 
lie in the $k$--linear space spanned by $f_{r+1}, \ldots, f_s$
as if this is the case the statement is trivial.  

Let $ y_0, \ldots, y_{n+r-s-1} \in S$ be polynomials of degree one 
which complete  $f_{r+1}, \ldots, f_s$ to a linear change of variables. 
We suppose in addition that $y_0=x_0$. Then the natural inclusion 
$k[y_0, \ldots, y_{n+r-s-1}] \hookrightarrow k[x_0,\ldots, x_n]/
(f_{r+1}, \ldots, f_s)$ is an isomorphism. 
Let  $v_i$ be an homogeneous polynomial in $k[y_0, \ldots, y_{n+r-s-1}]$
such that  $v_i \equiv  f_i $ $ \ \mod (f_{r+1}, \ldots, f_s)$ 
for $i=1, \ldots, r$. 
Then $x_0$ lies in the radical of the ideal $(v_1, \ldots, v_r)$ of 
$k[y_0, \ldots, y_{n+r-s-1}]$ and $\deg v_i \le d_i $ holds
for $i=1, \ldots, r$. 

Let $E$ denote the integer $2 \, \prod_{i=1}^r \, \deg v_i $ so that 
$E\le D:= 2 \, d_s\, \prod_{i=1}^{\min\{ n,s\}-1 }\, d_i $. 
Then $x_0^D\in (v_1, \ldots, v_r)$ from where it follows that 
$x_0^D \in (f_1, \ldots, f_s)$ as stated.
\end{proof}

Then Theorem 3 follows from this result by homogenizating the input 
polynomials and by considering the degree of the polynomials in a 
representation of $x_0^D $.

\vspace{3mm}

Now we are going to prove Theorem 4. First we recall the definition of 
algebraic degree of a polynomial system. 

Let $g_1, \ldots, g_s\in k[x_1, \ldots, x_n]$ be polynomials 
without common zeros in $\AA^n$. 
Given a $s\times s $ matrix $ \lambda=(\lambda_{ij})_{ij}  $
with entries in $\overline{k}$ we denote by 
$h_i(\lambda)$ the linear combinations $\sum_j \, 
\lambda_{ij} \, g_j $ induced by $\lambda$ for $i=1,  \ldots, s$. 
We consider the set of $s\times s $ matrices 
$\Gamma$ such that for any $\lambda$ in $\Gamma$ the polynomials 
$h_1(\lambda), \ldots, h_{t-1}(\lambda)$ form a
regular sequence in $\overline{k}[x_1, \ldots, x_n]$ and 
$1\in (h_1(\lambda), \ldots, h_t(\lambda))$ for some 
$t=t(\lambda) \le \min\{n,s\}$.  

For each $\lambda \in \Gamma$ and $i=1, \ldots, t-1$ we denote by  
$J_i(\lambda) \subseteq k[x_0, \ldots, x_n]$ the homogenization of 
the ideal $(h_1(\lambda ), \ldots, h_i(\lambda ))$. 
Then let $\delta(\lambda)$ denote the maximum degree 
of the homogeneous ideals $J_i(\lambda)$ for  $i=1 , \ldots, t-1$. 

The {\it algebraic degree} of the polynomial system $g_1, \ldots, g_s $ 
is defined as the minimum of $ \delta(\lambda)$ over all matrices  
$\lambda \in \Gamma$. 

The notion of geometric degree of \cite{KrSaSo96} and \cite{Sombra96} 
is defined in
an analogous way as the minimum of $\delta(\lambda)$ for $\lambda \in 
\Gamma$, with the additional hypothesis that the ideals 
$J_i(\lambda) $ are radical for $i=1, \ldots, t-1$. Another 
difference is that in the case when the characteristic of $k$ is positive
the polynomials $h_j (\lambda)$ are taken as linear combinations
of the polynomials $\{x_j\, f_i\}_{i j}$. 

The notion of geometric degree of $\cite{GiHeMoMoPa95}$ is similar to 
that of \cite{KrSaSo96}, \cite{Sombra96}, the only difference is 
that it is not
defined as a minimum but as the value of $\delta(\lambda)$ for a 
generic choice of $\lambda$. 

Thus the algebraic degree is bounded by the geometric degree, whichever
version of the later one we consider. 
The following example shows that in fact it can be much smaller. 
It is a variant of \cite[Example]{KrSaSo96}. 

\begin{exmp} \label{exmp3.1}
{\rm
Let us consider the polynomial system 
\[
f_1:=1-x_1x_2^d, \ f_2:=x_2-x_3^d, \ \ldots, \ f_{n-1}:=
x_{n-1}-x_n^d, \ f_n:=x_n^2  
\]
for some $d\ge 3$. 
It is easy to check that these polynomials have no common zero in $\AA^n$. 
Let  us denote by $\delta_g$ and by $\delta_a$ 
the geometric degree ---in  the sense of 
\cite{KrSaSo96}, \cite{Sombra96}--- and the algebraic degree 
of this polynomial system, respectively. 
We are going to compute both integers for this particular example. 

Let $\lambda$ be a $l\times n$ matrix with entries in $\overline{k}$
for some $l\le n$, and let $h_i := \sum_j \lambda_{ij} \, f_j$ be 
the induced linear combinations for $i=1, \ldots, n$. 

We apply to this matrix the Gauss elimination method by rows, by using as 
pivots the columns of $\lambda$ by ascending order. 
Invertible elementary operations by rows produces polynomials 
$q_1, \ldots, q_l$ which generate the same ideal than $h_1, \ldots, h_l$. 
For our particular polynomial system this corresponds to the 
successive elimination of the variables $x_1, \ldots, x_n$ in the 
equations $h_1, \ldots, h_l$. 
Thus in the case when $l\le n-1$ the affine variety defined by 
$h_1, \ldots, h_l$ can be parametrized 
by expressing the pivot variables as rational functions of the 
free ones. 
It follows that in this case the ideal $(h_1, \ldots, h_l)$ 
has dimension at least $n-l$. 

First we consider the geometric degree of this system. 
The polynomials $f_1, \ldots, f_n$ form a weak regular sequence, 
$1\in (f_1, \ldots, f_n)$ and the ideal $(f_1, \ldots, f_i)$ 
is radical for $i=1, \ldots, n-1$. 
We have that $\ \deg V(f_1, \ldots, f_i) = d^i \ $ for $i=1, \ldots, n-1$
from where it follows that $\delta_g \le d^{n-1}$. 

On the other hand, let $\lambda$ be a $l\times n $ matrix
with entries in $\overline{k}$ and let $h_1, \ldots, h_l$ be the corresponding 
linear combinations. 
Assume that  $h_1, \ldots, h_l$ is a 
weak regular sequence, $1\in (h_1, \ldots, h_l)$ and that 
$(h_1, \ldots, h_i) $ is a radical ideal for $i=1, \ldots, l$. 
In particular we get that $l$ equals $n$ and that 
$(h_1, \ldots, h_{n-1})$ is a one dimensional radical ideal. 

We apply the elimination method described above to the matrix 
formed by the first $n-1$ rows of $\lambda$ and we denote by 
$q_1, \ldots, q_{n-1}$ the obtained polynomials. 
We affirm that no columns $i$ fails to be a pivot for $i=1, \ldots, n-1$. 
Suppose that this is not the case. 
Then $q_1, \ldots, q_{n-2}$ are polynomials
which do not depend on $x_n$ and which span a one dimensional 
ideal in $\overline{k}[x_1, \ldots, x_{n-1}]$. 
In addition $q_{n-1}=x_n^2$ and so the ideal $(q_1, \ldots, q_{n-1})
\subseteq \overline{k}[x_1, \ldots, x_n]$ is not radical, contradicting 
our assumption. 

Thus there exist scalars $a_1, \ldots, a_{n-1} \in \overline{k}$ 
such that $\ q_i =f_i + a_i \, f_n \ $ for $i=1,\ldots, n-1$. 
We deduce that the variety $V(h_1, \ldots, h_{n-1})$ can
be parametrized by a map $\varphi: \AA^1  \to \AA^n$
defined by $ t\mapsto \varphi(t) = (\varphi_1(t), 
\ldots, \varphi_n(t))$, where $\varphi_i \in \overline{k}(t) $ 
is a rational function of degree $d^{n-i}$ for $i=1, \ldots, n$. 
We get that $ \ \deg V(h_1, \ldots, h_{n-1}) =d^{n-1} \ $ from where it 
follows the lower bound $\ \delta_g \ge d^{n-1}$. 
Combining this with the previous estimation we conclude that 
$\delta_g =d^{n-1}$. 

Now we consider the algebraic degree of the system. 
The polynomials $f_n, \ldots, f_1$ form a weak regular sequence
and $1\in (f_n, \ldots, f_1)$. 
We have that 
$(f_n, \ldots, f_{n-i+1})=( x_n^2, x_{n-1}, \ldots,x_{n-i+1})$  for 
$i=1,  \ldots,  n$ from where it follows that $\delta_a \le 2$. 
In addition, any nontrivial linear combination of $f_1, \ldots, 
f_n$ has degree at least two and so $\delta_a \ge 2$. 
We conclude that $\delta_a =2$. 
}
\end{exmp}

We obtain the following degree bound by direct 
application of Lemmas \ref{lem3.3} and \ref{lem3.4}. 

\begin{thm} \label{thm3.2}
Let be given homogeneous polynomials 
$f_1, \ldots, f_s \in k[x_0,\ldots, x_n]$ such that $x_0$ lies in the radical
of the ideal $(f_1, \ldots, f_s)$. Let $f_i^a$ denote the affinization 
of $f_i$ for $i=1, \ldots, s$. 
Let $d$ denote the maximum degree of $f_i$ for $i=1, \ldots, s$ and let 
$\delta $ denote the degree of the polynomial system $f_1^a, 
\ldots, f_s^a$. Then 
\[
x_0^D \in (f_1, \ldots, f_s)
\]
holds, with $ \ D:= \min\{ n, s\}^2 \, d \, \delta$. 
\end{thm}

\begin{proof}{}
\end{proof}

Then Theorem 4 follows from this result in the same way 
that Theorem 3 was derived from Theorem \ref{thm3.1}. 

If we apply this degree bound to the previous example we obtain
that there exist $g_1, \ldots, g_n\in k[x_1, \ldots, x_n]$ satisfying
\[
1=g_1f_1+ \cdots + g_nf_n,
\]
with $\ \deg g_if_i \le 2 \, n^2 \, d \ $ for $i=1, \ldots, s$. 
In fact we have the identity
\[
1= f_1 + x_1 \, x_2^{d-1} \, f_2+ x_1 \, x_2^{d-1}\, x_3^{d-1}\, f_3+  
	\cdots + x_1 \, x_2^{d-1} \cdots x_{n-1}^{d-1}\, x_n^{d-2}\,  f_n 
\]

\vspace{10mm}

%------------ Agradecimientos -----------------------------------------

\paragraph*{Acknowledgements.} This work originated in 
several conversations with Bernd Sturmfels. He motivated me to 
think about the sparse Nullstellensatz and suggested me several 
lines to approach it. Special thanks are due to him. I am
also grateful to Alicia Dickenstein and Joos Heintz for helpful 
discussions and suggestions, and to Pablo Solern\'o for providing
me a counterexample to a conjecture in an early version of this paper. 

I also thank the Departments of Mathematics of the 
Universities of Alcal\'a and of Cantabria, Spain, where part of 
this paper was written during a stay in the spring of 1997.

% ----------- Referencias -----------------------------------------------

\typeout{References}

\end{document}